\documentclass[
12pt]{article}
\usepackage{amssymb}
\usepackage{amsmath}
\usepackage{mathrsfs}
\usepackage{color}
\numberwithin{equation}{section}
\newtheorem{theo}{Theorem}[section]

\newtheorem{rem}[theo]{\it Remark}

\newtheorem{cor}[theo]{Corollary}
\newtheorem{lemma}[theo]{Lemma}

\def\pp{{\bf P}}
\def\ee{{\bf E}}
\def\dd{{\rm d}}

\def\LL{\mbox{\boldmath$Lin$}}
\def\Box{\rule{6pt}{6pt}}

\def\NR{\mathbb{R}}

\def\NZ{\mathbb{Z}}
\begin{document}
\title{A resolvent-type method for estimating time integrals of
quadratic fluctuations in weakly asymmetric exclusion}
\author{Sigurd Assing\\
Dept.\ of Statistics, The University of Warwick, Coventry CV4 7AL, UK\\
e-mail: s.assing@warwick.ac.uk}
\date{}
\maketitle
\begin{abstract}
For general right processes $\eta_s$ in a stationary state $\nu$,
under fairly weak conditions,
it is shown that 
$\int_0^T\hspace{-3pt}\dd t\,\ee_\nu[\int_0^{t}\hspace{-1pt}
V(s,\eta_{s\varepsilon^{-\kappa}})\dd s]^2
\le c_T\hspace{-2pt}\int_0^T\hspace{-3pt}
\left(\rule{0pt}{10pt}\right.\!
V(s,\cdot)
\;\rule[-4pt]{1pt}{14pt}\;
(1-\varepsilon^{-\kappa}L^{\rm sym})^{-1}V(s,\cdot)
\!\left.\rule{0pt}{10pt}\right)\!\dd s$
where $L^{\rm sym}$ denotes the symmetric part of the generator of $\eta_s$
on the Hilbert space $L^2(\nu)$ with inner product $(\cdot\,|\,\cdot)$.
Compared to
$\int_0^T\hspace{-3pt}
\left(\rule{0pt}{10pt}\right.\!
V(s,\cdot)
\;\rule[-4pt]{1pt}{14pt}\;
(-\varepsilon^{-\kappa}L^{\rm sym})^{-1}V(s,\cdot)
\!\left.\rule{0pt}{10pt}\right)\!\dd s$
which is often used in this context,
the advantage is that
$(1-\varepsilon^{-\kappa}L^{\rm sym})^{-1}V$
always exists for bounded measurable $V$.
As a consequence one obtains useful estimates of time integrals of
$\varepsilon$-scaled  quadratic fluctuations $V_\varepsilon$ build from 
$V_\#(\eta)=(\eta(0)-1/2)(\eta(1)-1/2)$ 
in the case of $\sqrt{\varepsilon}$-asymmetric exclusion
with $\nu$ being the symmetric Bernoulli product measure on
$\{0,1\}^\NZ$.
\end{abstract}
%
\noindent
{\large KEY WORDS}\hspace{0.4cm}right process, 
resolvent, generator, scaling limit, replacement lemma

\vspace{10pt}
\section{Motivation and Summary}
Assume that a particle system $\eta_s(x),\,x\in\NZ^d,\,s\ge 0$, 
gives raise to a scaled field of type 
\begin{equation}\label{field}
Y_s^\varepsilon\,=\,
\varepsilon^{\lambda_d}
\sum_x\frac{\eta_{s\varepsilon^{-\kappa}}(x)-a}{\chi}
\,\delta_{\varepsilon x-bs\varepsilon^{-\tilde{\kappa}}}
,\quad s\ge 0,
\end{equation}
where $\delta_{\varepsilon x-bs\varepsilon^{-\tilde{\kappa}}}$ 
denotes the Dirac measure
concentrated in the macroscopic point 
$\varepsilon x-bs\varepsilon^{-\tilde{\kappa}}$
and one wishes to understand the limiting behaviour, 
$\varepsilon\downarrow 0$,
of this field. 
Usually it follows from the martingale problem 
for the strong Markov process $\eta_s(x)$
that there is an approximate equation
for $Y_s^\varepsilon$, $\varepsilon$ small, which often reads like
\begin{equation}\label{approx equ}
\dd Y_s^\varepsilon(G)\,\sim\,
Y_s^\varepsilon({\cal A}G)\,\dd s\,+\,
V_\varepsilon^G(s,\xi_{s\varepsilon^{-\kappa}})\,\dd s\,+\,
\dd M_s^{G,\varepsilon}
\end{equation}
where $G$ is a smooth test function on $\NR^d$ with compact support,
${\cal A}$ a partial differential operator, 
$V_\varepsilon^G$ is for fixed $\varepsilon,G$ a bounded measurable function,
$\xi_{s\varepsilon^{-\kappa}}$ stands for
$(\eta_{s\varepsilon^{-\kappa}}(x)-a)/\chi$
and $M^{G,\varepsilon}$ denotes a martingale.
For further analysis of this approximate equation it becomes necessary
to express $V_\varepsilon^G(s,\xi_{s\varepsilon^{-\kappa}})$ in terms
of the field $Y_s^\varepsilon$, that is, one wants to replace
$$\int_0^t V_\varepsilon^G(s,\xi_{s\varepsilon^{-\kappa}})\,\dd s
\quad\mbox{by}\quad
\int_0^t F(s,Y_s^\varepsilon,G)\,\dd s$$
in some sense where the functional $F$ has to be found
(see \cite{KL1999} for a good review of this method). Typically the
difference
$$\int_0^t F(s,Y_s^\varepsilon,G)\,\dd s
\,-\int_0^t V_\varepsilon^G(s,\xi_{s\varepsilon^{-\kappa}})\,\dd s
\quad\mbox{simplifies to}\quad
\sum_{i=1}^m
\int_0^{t}
V_\varepsilon^{G,i}(s,\xi_{s\varepsilon^{-\kappa}})\,\dd s$$
and the task is to estimate an appropriate norm of
$t\mapsto\int_0^{t}\hspace{-1pt}
V_\varepsilon(s,\eta_{s\varepsilon^{-\kappa}})\dd s$
where $V_\varepsilon(s,\eta)$ stands for one of the functions
$V_\varepsilon^{G,i}(s,(\eta-a)/\chi),\,i=1,\dots,m$.
 
First, for an arbitrary but fixed $\beta>0$ 
and a finite time horizon $T$,
it follows from Lemma \ref{lem_inhomo}
in Section \ref{section resolv method} that
\begin{equation}\label{inhomo}
\int_0^T\hspace{-5pt}\dd t\,\ee_\nu\,[\int_0^{t}\hspace{-1pt}
V_\varepsilon(s,\eta_{s\varepsilon^{-\kappa}})\,\dd s\,]^2
\,\le\,
\frac{2e^{\beta T}}{\beta}\,\varepsilon^{2\kappa}
(\tilde{V}_\varepsilon\,|\,
G_{\frac{\beta}{2}\varepsilon^{\kappa}}\tilde{V}_\varepsilon)_{L^2(\dd s\otimes\dd\nu)}
\end{equation}
where $\tilde{V}_\varepsilon$ stands for the function
$(s,\eta)\mapsto e^{-\frac{\beta}{2}s\varepsilon^{\kappa}}
V_\varepsilon(s\varepsilon^{\kappa},\eta)$
and $(G_\alpha)_{\alpha>0}$ denotes the strongly continuous contraction resolvent
associated with the process $(s,\eta_s)_{s\ge 0}$
on the Hilbert space $H=L^2(\dd s\otimes\dd\nu)$ 
assuming that there exists an invariant state $\nu$ 
of the system $(\eta_s)_{s\ge 0}$.
Notice that (\ref{inhomo}) is valid in the context of general right
processes.

The observation is now that in many cases one has the inequality
\begin{equation}\label{not new}
(u\,|\,G_{\alpha}u)
\,\le\,
\left(\rule{0pt}{10pt}\right.\!
u
\;\rule[-4pt]{1pt}{14pt}\;
(\alpha-L_{s})^{-1}\,u
\!\left.\rule{0pt}{10pt}\right),
\quad u\in H,
\end{equation}
for all $\alpha>0$ by abstract theory on resolvents
where $L_{s}$ stands for the symmetric part of the generator of the
resolvent $(G_\alpha)_{\alpha>0}$ in $H$.

Second, choosing $\alpha=\frac{\beta}{2}\varepsilon^{\kappa}$ and
$u=\tilde{V}_\varepsilon$ in (\ref{not new})
and applying (\ref{inhomo}) yields
\begin{equation}\label{inhomo separated}
\int_0^T\hspace{-5pt}\dd t\,
\ee_\nu\,[\int_0^{t}\hspace{-1pt}
V_\varepsilon(s,\eta_{s\varepsilon^{-\kappa}})\,\dd s\,]^2
\,\le\,
\frac{2e^{\beta T}}{\beta}\,
\varepsilon^\kappa\hspace{-3pt}\int_0^T\hspace{-5pt}
\left(\rule{0pt}{10pt}\right.\!
V_\varepsilon(s,\cdot)
\;\rule[-4pt]{1pt}{14pt}\;
(\frac{\beta\varepsilon^\kappa}{2}-L^{\rm sym})^{-1}V_\varepsilon(s,\cdot)
\!\left.\rule{0pt}{10pt}\right)_{L^2(\nu)}\dd s
\end{equation}
where $L^{\rm sym}$ denotes the symmetric part of the generator of the
process $(\eta_s)_{s\ge 0}$ in $L^2(\nu)$.
The details of how to replace $L_s$ by $L^{\rm sym}$ 
are explained by Lemma \ref{decouple} 
in Section \ref{section resolv method}. 
 
Remark that applying 
Kipnis-Varadhan's inequality\footnote{See \cite[Lemma 4.3]{CLO2001}
for the version used here.} which is widely used in the context of
particle systems would yield
$$\int_0^T\hspace{-5pt}\dd t\,
\ee_\nu\,[\int_0^{t}\hspace{-1pt}
V_\varepsilon(s,\eta_{s\varepsilon^{-\kappa}})\,\dd s\,]^2
\,\le\,
14T\,\varepsilon^\kappa\hspace{-3pt}\int_0^T\hspace{-5pt}
\left(\rule{0pt}{10pt}\right.\!
V_\varepsilon(s,\cdot)
\;\rule[-4pt]{1pt}{14pt}\;
(-L^{\rm sym})^{-1}V_\varepsilon(s,\cdot)
\!\left.\rule{0pt}{10pt}\right)_{L^2(\nu)}\dd s$$
instead of (\ref{inhomo separated}).
But, in the important case where
$(\eta_s)_{s\ge 0}$ is a simple one-dimensional exclusion process and
$\nu$ is the symmetric Bernoulli product measure on $\{0,1\}^\NZ$,
the right-hand side of the last inequality 
is infinite for $V_\#(\eta)=(\eta(0)-1/2)(\eta(1)-1/2)$
which is the quadratic part of the normalised current and the most
basic quadratic fluctuation.
So Kipnis-Varadhan's inequality cannot be used 
to estimate time integrals of
$\varepsilon$-scaled  quadratic fluctuations 
$V_\varepsilon$ build from $V_\#$.

However, the right-hand side of the inequality (\ref{inhomo separated})
is always finite for bounded measurable functions $V_\varepsilon$.
Hence, when substituting $V_\varepsilon^{G,i}$ for
$V_\varepsilon,\,i=1,\dots,m$, this inequality 
gives a tool for how to show
\begin{equation}\label{tool}
\lim_{\varepsilon\to 0}\;
\int_0^T\hspace{-3pt}\dd t\,
\ee_\nu\left(
\int_0^t F(s,Y_s^\varepsilon,G)\,\dd s
\,-\,
\int_0^t V_\varepsilon^G(s,\xi_{s\varepsilon^{-\kappa}})\,\dd s
\right)^{\hspace{-2pt}2}\,=\,0.
\end{equation}
Remark that replacing
$\int_0^t V_\varepsilon^G(s,\xi_{s\varepsilon^{-\kappa}})\,\dd s$
by 
$\int_0^t F(s,Y_s^\varepsilon,G)\,\dd s$
in the sense of the above limit
is rather \underline{weak} since the replacement does \underline{not} 
hold for every $t\in[0,T]$ but only for an average over $t\in[0,T]$.
Nevertheless, as recently shown in \cite{A2012}, 
this weak form of a replacement 
is still sufficient for deriving a meaningful
equation which could be used as the limit of (\ref{approx equ}).
 
Section \ref{section resolv method} presents the resolvent method
which is based on 
(\ref{inhomo}),(\ref{not new}),(\ref{inhomo separated}).
A detailed proof of the inequality (\ref{not new})
is given in a general setting (see Corollary \ref{my inequ}).
The result as such cannot be new. However, the author could not
find a reference. The proof is based on a variational formula
(see Lemma \ref{variational}) which was also used in the proof of
\cite[Lemma 2.1]{LQSY2004} but without explicit proof and only in
the framework of exclusion processes.
The detailed proof is added for completeness and for
having a good account on the precise conditions needed.

In Section \ref{section application} the resolvent method is applied 
in the case of $\sqrt{\varepsilon}$-asymmetric one-dimensional simple
exclusion in equilibrium. The field (\ref{field}) of interest is the diffusively
scaled density fluctuation field. In the corresponding equation
(\ref{approx equ}), ${\cal A}$ is the one-dimensional Laplacian and
the bounded functions $V_\varepsilon^G$ are 
$\varepsilon$-scaled  quadratic fluctuations 
build from $V_\#$ introduced above.
The key result is Lemma \ref{key result} which leads to
a replacement of the time integrals of these quadratic fluctuations in
the sense of (\ref{tool}), see Corollary \ref{weak replacement}.
Remark that the density fluctuations in
$\sqrt{\varepsilon}$-asymmetric exclusion are related to non-trivial
distributions like the Tracy-Widom distribution hence they are a good
`medium' for testing techniques. Proving Lemma \ref{key result} using a
resolvent-type method can be considered to be such a test.\\

\noindent
{\it Acknowledgement}. 
The author thanks Wilhelm Stannat for helpful discussions.
\section{A Resolvent Method}\label{section resolv method}
Let $X$ be a Hausdorff topological space and assume that the
Borel-$\sigma$-algebra on $X$ is equal to the $\sigma$-algebra
generated by the set of all continuous functions on $X$.

Denote by 
$(\Omega,{\cal F},\pp_{\!\!\eta},\eta\in X,(\eta_s)_{s\ge 0})$ 
a right process with state space $X$, infinite life time and
corresponding filtration ${\cal F}_s,\,s\ge 0$, satisfying the usual
conditions (see \cite{S1988} for a good account on general right
processes). Assume that there exists a measure $\nu$ on $X$ which is
an invariant state of this right process and denote by $\ee_\nu$ the
expectation operator given by the probability measure
$\int_X\pp_{\!\!\eta}\,\nu(\dd\eta)$.

Obviously, the pair $(s,\eta_s)_{s\ge 0}$ gives another right process
$(\Omega,{\cal F},\pp_{\!\!s,\eta},
(s,\eta)\in[0,\infty)\times X,(s,\eta_s)_{s\ge 0})$
corresponding to the same filtration ${\cal F}_s,\,s\ge 0$, such that
$$\pp_{\!\!s,\eta}\left(\rule{0pt}{10pt}\right.
\{\eta_s=\eta\}
\left.\rule{0pt}{10pt}\right)\,=\,1
\quad\mbox{for all}\quad(s,\eta)\in[0,\infty)\times X.$$
Define the transition semigroup and the resolvent of $(s,\eta_s)_{s\ge 0}$
by
$$p_r V(s,\eta)\,=\,\ee_{s,\eta}V(s+r,\eta_{s+r}),\quad r\ge 0,$$
and
$$R_\alpha V(s,\eta)\,=\,
\int_0^\infty\dd r\,e^{-\alpha r}\,p_r V(s,\eta),\quad\alpha>0,$$
respectively, where $V$ is an arbitary bounded measurable function on
$[0,\infty)\times X$.

Denote by $\ell$ the Lebesgue measure on $[0,\infty)$ and notice that
$\ell\otimes\nu$ is an excessive measure on $[0,\infty)\times X$
with respect to
$(\Omega,{\cal F},\pp_{\!\!s,\eta},
(s,\eta)\in[0,\infty)\times X,(s,\eta_s)_{s\ge 0})$
because
$$\int p_r V\,\dd(\ell\otimes\nu)\,\le\,
\int V\,\dd(\ell\otimes\nu),\quad r\ge 0,$$
for all non-negative measurable functions 
$V$ on $[0,\infty)\times X$.
As a consequence, see Section IV.2 in \cite{MR1992} for the details,
there exists a strongly continuous contraction resolvent
$(G_\alpha)_{\alpha>0}$ on $L^2(\ell\otimes\nu)$ such that
$G_\alpha V$ is an $(\ell\otimes\nu)$-version of $R_\alpha V$ for all $\alpha>0$
and all bounded functions $V$ in $L^2(\ell\otimes\nu)$. 
\begin{lemma}\label{lem_inhomo}
Fix $\beta>0$, let $T>0$ be a finite time horizon and consider the
process $(\eta_{cs})_{s\ge 0}$ time-scaled by a factor $c>0$. Then
$$\int_0^T\hspace{-5pt}\dd t\,\ee_\nu\,[\int_0^{t}\hspace{-1pt}
V(s,\eta_{cs})\,\dd s\,]^2
\,\le\,
\frac{2e^{\beta T}}{\beta c^2}\,
(\tilde{V}\,|\,
G_{\!\frac{\beta}{2c}}\tilde{V})_{L^2(\ell\otimes\nu)}$$
for all bounded functions $V$ in $L^2(\ell\otimes\nu)$
where $\tilde{V}(s,\eta)=e^{-\frac{\beta}{2}s/c}\,V(s/c\,,\eta)$.
\end{lemma}
{\it Proof}. Choose a bounded function $V$ in $L^2(\ell\otimes\nu)$
and set $\hat{V}(s,\eta)=V(s/c\,,\eta)$. Then:
$$\int_0^T\hspace{-5pt}\dd t\,\ee_\nu\,[\int_0^{t}\hspace{-1pt}
V(s,\eta_{cs})\,\dd s\,]^2
\,\le\,
e^{\beta T}
\int_0^\infty\hspace{-5pt}\dd t\,e^{-\beta t}\,\ee_\nu\,[\int_0^{t}\hspace{-1pt}
\hat{V}(cs,\eta_{cs})\,\dd s\,]^2$$
$$=\,\frac{2e^{\beta T}}{c^2}
\int_0^\infty\hspace{-5pt}\dd t\,e^{-\beta t}\,\ee_\nu
\int_0^{ct}\hspace{-1pt}\dd s\,\hat{V}(s,\eta_{s})
\int_s^{ct}\hspace{-1pt}\dd r\,\hat{V}(r,\eta_{r})$$
$$=\,\frac{2e^{\beta T}}{c^2}
\int_0^\infty\hspace{-5pt}\dd t\,e^{-\beta t}\hspace{-5pt}
\int_0^{ct}\hspace{-5pt}\dd s\int_s^{ct}\hspace{-5pt}\dd r\,
\ee_\nu\hat{V}(s,\eta_{s})\,p_{r-s}\hat{V}(s,\eta_{s})$$
$$=\,\frac{2e^{\beta T}}{c^2}
\int_0^\infty\hspace{-5pt}\dd t\,e^{-\beta t}\hspace{-5pt}
\int_0^{ct}\hspace{-5pt}\dd s\int_0^{ct-s}\hspace{-10pt}\dd r
\int\dd\nu\,\hat{V}(s,\cdot)\,p_{r}\hat{V}(s,\cdot)$$
\begin{equation}\label{no t}
=\,\frac{2e^{\beta T}}{\beta c^2}
\int_0^{\infty}\hspace{-5pt}\dd s\int_0^{\infty}\hspace{-5pt}\dd r
\,e^{-\beta(s+r)/c}\hspace{-5pt}
\int\dd\nu\,\hat{V}(s,\cdot)\,p_{r}\hat{V}(s,\cdot)
\end{equation}
$$=\,\frac{2e^{\beta T}}{\beta c^2}
\int_0^{\infty}\hspace{-5pt}\dd s\int_0^{\infty}\hspace{-5pt}\dd r
\,e^{-\beta(s+r)/c}\,
\ee_\nu\,\hat{V}(s,\eta_{s})\hat{V}(s+r,\eta_{s+r})$$
$$=\,\frac{2e^{\beta T}}{\beta c^2}
\int_0^{\infty}\hspace{-5pt}\dd s\int_0^{\infty}\hspace{-5pt}\dd r
\,e^{-\frac{\beta}{2c}r}\,
\ee_\nu\,\tilde{V}(s,\eta_{s})\tilde{V}(s+r,\eta_{s+r})$$
$$=\,\frac{2e^{\beta T}}{\beta c^2}
\int_0^{\infty}\hspace{-5pt}\dd s
\int\dd\nu\,\,\tilde{V}(s,\cdot)\hspace{-3pt}
\int_0^{\infty}\hspace{-5pt}\dd r\,e^{-\frac{\beta}{2c}r}\,
p_{r}\tilde{V}(s,\cdot)
\,=\,
\frac{2e^{\beta T}}{\beta c^2}\,
(\tilde{V}\,|\,
G_{\!\frac{\beta}{2c}}\tilde{V})_{L^2(\ell\otimes\nu)}
\eqno\Box$$
\begin{rem}\rm\label{no t remark}
\begin{itemize}\item[(i)]
In the case were $V$ is time independent
it easily follows from (\ref{no t}) that
$$\int_0^T\hspace{-5pt}\dd t\,\ee_\nu\,[\int_0^{t}\hspace{-1pt}
V(\eta_{cs})\,\dd s\,]^2
\,\le\,
\frac{2e^{\beta T}}{\beta^2 c}\,
({V}\,|\,
G_{\!\beta/c}{V})_{L^2(\nu)}$$
because $V=\hat{V}$ and 
$\int_0^{\infty}\hspace{-2pt}\dd s\,e^{-\beta s/c}
=c/\beta$.
\item[(ii)]
Notice that one cannot apply (i) to the process $(s,\eta_s)_{s\ge 0}$
because $\ell\otimes\nu$ is not an invariant measure for this process.
\end{itemize}
\end{rem}

The remaining part of this section deals with the problem of
estimating the right-hand side of the inequality 
in Lemma \ref{lem_inhomo} by something which 
is more likely to be computable in an explicit way.

Let $(G_\alpha)_{\alpha>0}$ be a strongly continuous contraction
resolvent on a real Hilbert space $H$ with inner product
$(\cdot\,|\,\cdot)$. If $G_\alpha^\ast$ denotes the adjoint of $G_\alpha$ 
then $(G_\alpha^\ast)_{\alpha>0}$ is also a strongly continuous
contraction resolvent on $H$.
Denote by $(L,D(L))$ and $(L^\ast,D(L^\ast))$ the generators of
$(G_\alpha)_{\alpha>0}$ and $(G_\alpha^\ast)_{\alpha>0}$, respectively.
\begin{rem}\rm
If the co-generator $(L^\ast,D(L^\ast))$ is not the adjoint of $(L,D(L))$ 
then $L^\ast$ coincides with the adjoint of $(L,D(L))$ on 
$D(L^\ast)$, at least.
\end{rem}
The following assumption will play a
crucial role in the proof of Lemma \ref{variational} below.
$$\mbox{There exists $D_0\subseteq H$ which is a core for both
$(L,D(L))$ and $(L^\ast,D(L^\ast))$.}\eqno(A1)$$
Remark that there are unbounded operators $(L,D(L))$ 
on Hilbert spaces whose adjoints $(L^\ast,D(L^\ast))$ satisfy
$D(L)\cap D(L^\ast)=\{0\}$ in the worst case
hence there is something to check for this assumption to hold.

Assuming $(A1)$, the symmetric and antisymmetric parts of $L$ given by
$$L_s\,=\,\frac{L+L^\ast}{2}
\quad\mbox{and}\quad
L_a\,=\,\frac{L-L^\ast}{2}
\quad\mbox{respectively}$$
are defined on $D_0$. Of course $D_0$ is dense in $H$ since it is a core.
As $(L,D(L))$ and $(L^\ast,D(L^\ast))$ are both negative definite, 
$(L_s,D_0)$ is a symmetric, negative definite, densely defined operator. 
Hence it can be extended (Friedrich's extension for example) 
to a self-adjoint negative definite operator $(L_s,D(L_s))$ on $H$.
Every such extension generates a strongly continuous
contraction resolvent $((\alpha-L_s)^{-1})_{\alpha\ge 0}$
of self-adjoint positive definite operators on $H$.
\begin{lemma}\label{variational}
Assume $(A1)$ and
fix both $\alpha>0$ as well as an arbitrary self-adjoint extension
of $(L_s,D_0)$. Then
$$(u\,|\,G_\alpha u)\,=\,\sup_{v\in D_0}
\left\{\rule{0pt}{12pt}\right.
2\,(u\,|\,v)-
\left(\rule{0pt}{10pt}\right.\!
(\alpha-L)v
\;\rule[-4pt]{1pt}{14pt}\;
(\alpha-L_s)^{-1}(\alpha-L)v
\!\left.\rule{0pt}{10pt}\right)
\left.\rule{0pt}{12pt}\right\}$$
for all $u\in H$.
\end{lemma}
{\it Proof}. Fix $u\in H$. As $D_0$ is a core for both, $(L,D(L))$ and
$(L^\ast,D(L^\ast))$, one can choose sequences $(v_n)_{n=1}^\infty$
and $(v_n^\ast)_{n=1}^\infty$ in $D_0$ such that
$$v_n\to\frac{G_\alpha}{2}u,\quad
(\alpha-L)v_n=:\frac{u_n}{2}\to\frac{u}{2},\quad
v_n^\ast\to\frac{G_\alpha^\ast}{2}u,\quad
(\alpha-L^\ast)v_n^\ast=:\frac{u_n^\ast}{2}\to\frac{u}{2},\quad$$
in $H$ when $n\to\infty$. Then
$$\left(\rule{0pt}{10pt}\right.
(\alpha-L)(v_n+v_n^\ast)
\;\rule[-4pt]{1pt}{14pt}\;
(\alpha-L_s)^{-1}(\alpha-L)(v_n+v_n^\ast)
\left.\rule{0pt}{10pt}\right)
\longrightarrow
(u\,|\,G_\alpha u),\quad n\to\infty.$$
Indeed
\begin{eqnarray*}
&&
\left(\rule{0pt}{10pt}\right.
(\alpha-L)(v_n+v_n^\ast)
\;\rule[-4pt]{1pt}{14pt}\;
(\alpha-L_s)^{-1}(\alpha-L)(v_n+v_n^\ast)
\left.\rule{0pt}{10pt}\right)\\
&=&
\left(\rule{0pt}{10pt}\right.
\frac{\alpha-L}{2}\,(v_n+v_n^\ast)
\;\rule[-4pt]{1pt}{14pt}\;
(\alpha-L_s)^{-1}u_n
\left.\rule{0pt}{10pt}\right)
+\left(\rule{0pt}{10pt}\right.
(\alpha-L)(v_n+v_n^\ast)
\;\rule[-4pt]{1pt}{14pt}\;
(\alpha-L_s)^{-1}\,\frac{\alpha-L}{2}\,2v_n^\ast
\left.\rule{0pt}{10pt}\right)
\end{eqnarray*}
where
\begin{eqnarray*}
&&
\left(\rule{0pt}{10pt}\right.
(\alpha-L)(v_n+v_n^\ast)
\;\rule[-4pt]{1pt}{14pt}\;
(\alpha-L_s)^{-1}\,\frac{\alpha-L}{2}\,2v_n^\ast
\left.\rule{0pt}{10pt}\right)\\
&=&
\left(\rule{0pt}{10pt}\right.
(\alpha-L)(v_n+v_n^\ast)
\;\rule[-4pt]{1pt}{14pt}\;
2v_n^\ast
\left.\rule{0pt}{10pt}\right)
-\left(\rule{0pt}{10pt}\right.
(\alpha-L)(v_n+v_n^\ast)
\;\rule[-4pt]{1pt}{14pt}\;
(\alpha-L_s)^{-1}\,\frac{\alpha-L^\ast}{2}\,2v_n^\ast
\left.\rule{0pt}{10pt}\right)\\
&=&
\left(\rule{0pt}{10pt}\right.
(v_n+v_n^\ast)
\;\rule[-4pt]{1pt}{14pt}\;
u_n^\ast
\left.\rule{0pt}{10pt}\right)
-\left(\rule{0pt}{10pt}\right.
\frac{\alpha-L}{2}\,(v_n+v_n^\ast)
\;\rule[-4pt]{1pt}{14pt}\;
(\alpha-L_s)^{-1}u_n^\ast
\left.\rule{0pt}{10pt}\right)
\end{eqnarray*}
hence
\begin{eqnarray*}
&&
\left(\rule{0pt}{10pt}\right.
(\alpha-L)(v_n+v_n^\ast)
\;\rule[-4pt]{1pt}{14pt}\;
(\alpha-L_s)^{-1}(\alpha-L)(v_n+v_n^\ast)
\left.\rule{0pt}{10pt}\right)\\
&=&
\left(\rule{0pt}{10pt}\right.
(v_n+v_n^\ast)
\;\rule[-4pt]{1pt}{14pt}\;
u_n^\ast
\left.\rule{0pt}{10pt}\right)
+\left(\rule{0pt}{10pt}\right.
\frac{\alpha-L}{2}\,(v_n+v_n^\ast)
\;\rule[-4pt]{1pt}{14pt}\;
(\alpha-L_s)^{-1}(u_n-u_n^\ast)
\left.\rule{0pt}{10pt}\right).
\end{eqnarray*}
Of course
$$\left(\rule{0pt}{10pt}\right.
(v_n+v_n^\ast)
\;\rule[-4pt]{1pt}{14pt}\;
u_n^\ast
\left.\rule{0pt}{10pt}\right)
\,=\,
\left(\rule{0pt}{10pt}\right.
u_n^\ast
\;\rule[-4pt]{1pt}{14pt}\;
(v_n+v_n^\ast)
\left.\rule{0pt}{10pt}\right)
\longrightarrow(u\,|\,G_\alpha u),\quad n\to\infty,$$
whereas
\begin{eqnarray*}
\left(\rule{0pt}{10pt}\right.
\frac{\alpha-L}{2}\,(v_n+v_n^\ast)
\;\rule[-4pt]{1pt}{14pt}\;
(\alpha-L_s)^{-1}(u_n-u_n^\ast)
\left.\rule{0pt}{10pt}\right)
&=&
\left(\rule{0pt}{10pt}\right.
\frac{\alpha-L}{2}\,v_n
\;\rule[-4pt]{1pt}{14pt}\;
(\alpha-L_s)^{-1}(u_n-u_n^\ast)
\left.\rule{0pt}{10pt}\right)\\
&-&
\left(\rule{0pt}{10pt}\right.
\frac{\alpha-L^\ast}{2}\,v_n^\ast
\;\rule[-4pt]{1pt}{14pt}\;
(\alpha-L_s)^{-1}(u_n-u_n^\ast)
\left.\rule{0pt}{10pt}\right)\\
&+&
\left(\rule{0pt}{10pt}\right.
(\alpha-L_s)v_n^\ast
\;\rule[-4pt]{1pt}{14pt}\;
(\alpha-L_s)^{-1}(u_n-u_n^\ast)
\left.\rule{0pt}{10pt}\right)
\end{eqnarray*}
converges to zero when $n\to\infty$ which is obvious for the first two
terms on the right-hand side and follows from
$$\left(\rule{0pt}{10pt}\right.
(\alpha-L_s)v_n^\ast
\;\rule[-4pt]{1pt}{14pt}\;
(\alpha-L_s)^{-1}(u_n-u_n^\ast)
\left.\rule{0pt}{10pt}\right)
\,=\,
\left(\rule{0pt}{10pt}\right.
v_n^\ast
\;\rule[-4pt]{1pt}{14pt}\;
(u_n-u_n^\ast)
\left.\rule{0pt}{10pt}\right)$$
for the last term. 

Altogether one obtains that
\begin{eqnarray*}
&&(u\,|\,G_\alpha u)\\
&=&
\lim_{n\to\infty}
\left\{\rule{0pt}{12pt}\right.
\hspace{-2pt}
2\left(\rule{0pt}{10pt}\right.
u
\;\rule[-4pt]{1pt}{14pt}\;
(v_n+v_n^\ast)
\left.\rule{0pt}{10pt}\right)
-\left(\rule{0pt}{10pt}\right.
(\alpha-L)(v_n+v_n^\ast)
\;\rule[-4pt]{1pt}{14pt}\;
(\alpha-L_s)^{-1}(\alpha-L)(v_n+v_n^\ast)
\left.\rule{0pt}{10pt}\right)
\hspace{-2pt}
\left.\rule{0pt}{12pt}\right\}\\
&\le&
\sup_{v\in D_0}
\left\{\rule{0pt}{12pt}\right.
2\,(u\,|\,v)-
\left(\rule{0pt}{10pt}\right.\!
(\alpha-L)v
\;\rule[-4pt]{1pt}{14pt}\;
(\alpha-L_s)^{-1}(\alpha-L)v
\!\left.\rule{0pt}{10pt}\right)
\left.\rule{0pt}{12pt}\right\}
\end{eqnarray*}
and it remains to show that
\begin{equation}\label{other direction}
\sup_{v\in D_0}
\left\{\rule{0pt}{12pt}\right.
2\,(u\,|\,v)-
\left(\rule{0pt}{10pt}\right.\!
(\alpha-L)v
\;\rule[-4pt]{1pt}{14pt}\;
(\alpha-L_s)^{-1}(\alpha-L)v
\!\left.\rule{0pt}{10pt}\right)
\left.\rule{0pt}{12pt}\right\}
\,\le\,
(u\,|\,G_\alpha u).
\end{equation}
Now choose $(v_n^\ast)_{n=1}^\infty\subseteq D_0$ such that
$$v_n^\ast\to G_\alpha^\ast u
\quad\mbox{and}\quad
(\alpha-L^\ast)v_n^\ast\to u
\quad\mbox{in $H$ when $n\to\infty$}$$
and remark that
$(\alpha-L_s)^{1/2}$ is well-defined on $D_0\subseteq D(L_s)$.
Then for every $v\in D_0$
\begin{eqnarray*}
(u\,|\,v)&=&
\left(\rule{0pt}{10pt}\right.
u
\;\rule[-4pt]{1pt}{14pt}\;
G_\alpha(\alpha-L)v
\left.\rule{0pt}{10pt}\right)
\,=\,
\left(\rule{0pt}{10pt}\right.
G_\alpha^\ast u
\;\rule[-4pt]{1pt}{14pt}\;
(\alpha-L)v
\left.\rule{0pt}{10pt}\right)\\
&=&
\lim_{n\to\infty}
\left(\rule{0pt}{10pt}\right.
v_n^\ast
\;\rule[-4pt]{1pt}{14pt}\;
(\alpha-L_s)^{1/2}(\alpha-L_s)^{-1/2}(\alpha-L)v
\left.\rule{0pt}{10pt}\right)\\
&=&
\lim_{n\to\infty}
\left(\rule{0pt}{10pt}\right.
(\alpha-L_s)^{1/2}v_n^\ast
\;\rule[-4pt]{1pt}{14pt}\;
(\alpha-L_s)^{-1/2}(\alpha-L)v
\left.\rule{0pt}{10pt}\right)\\
&\le&
\lim_{n\to\infty}
\sqrt{
\left(\rule{0pt}{10pt}\right.
v_n^\ast
\;\rule[-4pt]{1pt}{14pt}\;
(\alpha-L_s)v_n^\ast
\left.\rule{0pt}{10pt}\right)
\cdot
\left(\rule{0pt}{10pt}\right.
(\alpha-L_s)v
\;\rule[-4pt]{1pt}{14pt}\;
(\alpha-L_s)^{-1}(\alpha-L)v
\left.\rule{0pt}{10pt}\right)
}\\
&\le&
\lim_{n\to\infty}
\left[\rule{0pt}{11pt}\right.
\left(\rule{0pt}{10pt}\right.
v_n^\ast
\;\rule[-4pt]{1pt}{14pt}\;
(\alpha-L_s)v_n^\ast
\left.\rule{0pt}{10pt}\right)
+
\left(\rule{0pt}{10pt}\right.
(\alpha-L_s)v
\;\rule[-4pt]{1pt}{14pt}\;
(\alpha-L_s)^{-1}(\alpha-L)v
\left.\rule{0pt}{10pt}\right)
\left.\rule{0pt}{11pt}\right]/2\\
&=&
\left[\rule{0pt}{11pt}\right.
\left(\rule{0pt}{10pt}\right.
(u
\;\rule[-4pt]{1pt}{14pt}\;
G_\alpha u
\left.\rule{0pt}{10pt}\right)
+
\left(\rule{0pt}{10pt}\right.
(\alpha-L_s)v
\;\rule[-4pt]{1pt}{14pt}\;
(\alpha-L_s)^{-1}(\alpha-L)v
\left.\rule{0pt}{10pt}\right)
\left.\rule{0pt}{11pt}\right]/2
\end{eqnarray*}
since
$$(v_n^\ast\,|\,(\alpha-L_s)v_n^\ast)
=(v_n^\ast\,|\,(\alpha-L^\ast)v_n^\ast)
\longrightarrow
(G_\alpha^\ast u\,|\,u)
\,=\,
(u\,|\,G_\alpha u),\quad n\to\infty,$$
and (\ref{other direction}) follows.\hfill$\Box$
\begin{rem}\rm
For the proof of (\ref{other direction}), $D_0$ only needs
to be a core for $(L^\ast,D(L^\ast))$ but not for $(L,D(L))$.
Of course $D_0\subseteq D(L)$ must still be assumed.
\end{rem}
\begin{cor}\label{my inequ}
If there exists $D_0\subseteq H$ 
which is a core for $(L,D(L)),\,(L^\ast,D(L^\ast))$ and $(L_s,D(L_s))$ then
$$(u\,|\,G_\alpha u)
\,\le\,
\left(\rule{0pt}{10pt}\right.\!
u
\;\rule[-4pt]{1pt}{14pt}\;
(\alpha-L_s)^{-1}\,u
\!\left.\rule{0pt}{10pt}\right),\quad u\in H,$$
for every $\alpha>0$.
\end{cor}
{\it Proof}.
If $D_0$ is a core for $(L,D(L)),\,(L^\ast,D(L^\ast))$ and
$(L_s,D(L_s))$ then Lemma \ref{variational} implies
$$(u\,|\,G_\alpha u)\,\le\,\sup_{v\in D_0}
\left\{\rule{0pt}{12pt}\right.
2\,(u\,|\,v)-
\left(\rule{0pt}{10pt}\right.\!
v
\;\rule[-4pt]{1pt}{14pt}\;
(\alpha-L_s)v
\!\left.\rule{0pt}{10pt}\right)
\left.\rule{0pt}{12pt}\right\}
\,=\,
\left(\rule{0pt}{10pt}\right.\!
u
\;\rule[-4pt]{1pt}{14pt}\;
(\alpha-L_s)^{-1}\,u
\!\left.\rule{0pt}{10pt}\right).$$
Indeed
\begin{eqnarray*}
\left(\rule{0pt}{10pt}\right.\!
(\alpha-L)v
\;\rule[-4pt]{1pt}{14pt}\;
(\alpha-L_s)^{-1}(\alpha-L)v
\!\left.\rule{0pt}{10pt}\right)
&=&
\left(\rule{0pt}{10pt}\right.\!
v
\;\rule[-4pt]{1pt}{14pt}\;
(\alpha-L_s)v
\!\left.\rule{0pt}{10pt}\right)
+
\left(\rule{0pt}{10pt}\right.\!
L_a v
\;\rule[-4pt]{1pt}{14pt}\;
(\alpha-L_s)^{-1}L_a v
\!\left.\rule{0pt}{10pt}\right)\\
&\ge&
\left(\rule{0pt}{10pt}\right.\!
v
\;\rule[-4pt]{1pt}{14pt}\;
(\alpha-L_s)v
\!\left.\rule{0pt}{10pt}\right)
\end{eqnarray*}
for $v\in D_0$ since $(\alpha-L_s)^{-1}$ is positive definite
and finally one applies Lemma \ref{variational} in the special case
$L=L_s$.\hfill$\Box$
\begin{rem}\rm
Not every perturbation $(L,D_0)$ of a symmetric operator 
$(S,D_0)$ has $S$ as its symmetric part on $D_0$. 
An easy example demonstrating this fact
will be given in the next section, see Remark \ref{not sym part}.
So one has to be careful when checking the assumptions of 
Corollary \ref{my inequ}.
\end{rem}

Finally the impact of Corollary \ref{my inequ} on the situation
described in Lemma \ref{lem_inhomo} is discussed.
So $H=L^2(\ell\otimes\nu)$, $(G_\alpha)_{\alpha>0}$
is the strongly continuous contraction resolvent
associated with the right process $(s,\eta_s)_{s\ge 0}$ on $H$
and $(L,D(L))$ is the generator of $(G_\alpha)_{\alpha>0}$.
Denote by $(L^\eta,D(L^\eta))$ the generator of the strongly continuous
contraction resolvent associated with the right process 
$(\eta_s)_{s\ge 0}$ on $L^2(\nu)$. 
Corollary \ref{my inequ} suggests to develop the inequality 
given in Lemma \ref{lem_inhomo} by using
$$(\tilde{V}\,|\,
G_{\!\frac{\beta}{2c}}\tilde{V})_{L^2(\ell\otimes\nu)}
\,\le\,
(\tilde{V}\,|\,
(\frac{\beta}{2c}-L_s)^{-1}\tilde{V})_{L^2(\ell\otimes\nu)}$$
and one wants to simplify 
$(\beta/(2c)-L_s)^{-1}\tilde{V}$
in this specific situation.
Formally $L=\frac{\partial}{\partial s}+L^\eta$ and
$L^\ast=-\delta_0(s)-\frac{\partial}{\partial s}+(L^\eta)^\ast$
where $\delta_0$ denotes Dirac's delta function.
So, under certain conditions on $\tilde{V}$, one should have an
equality of the type
$$
[(\alpha-L_s)^{-1}\tilde{V}](s,\eta)
\,=\,
[(\alpha-L^{\rm sym})^{-1}\tilde{V}(s,\cdot)](\eta)
\,\stackrel{\mbox{\tiny def}}{=}\,g(s,\eta)
$$
where $L^{\rm sym}$ denotes 
the symmetric part of $L^\eta$ in $L^2(\nu)$.
The following lemma, first, lists sufficient conditions 
to ensure this and, second, presents the final
bound on the left-hand side of the inequality in Lemma \ref{lem_inhomo}.
\begin{lemma}\label{decouple}
\begin{itemize}\item[(i)] 
Fix $\alpha>0$.
If $g,L^\eta g,(L^\eta)^\ast g,\frac{\partial}{\partial s}g
\in L^2(\ell\otimes\nu)$,
$g(s,\cdot)\in D(L^\eta)\cap D((L^\eta)^\ast)$, $s\ge 0$,
and
$g(\cdot,\eta)\in C_0^1([0,\infty)),\,\eta\in X$,
then
$(\alpha-L_s)^{-1}\tilde{V}\,=\,g$.

\item[(ii)] Fix $\beta>0,\,T>0,\,c>0$ and assume that there exists
$D_0\subseteq L^2(\nu)$ which is a core for 
$(L^\eta,D(L^\eta))$, $((L^\eta)^\ast,D((L^\eta)^\ast))$
and $(L^{\rm sym},D(L^{\rm sym}))$. Then:
$$\int_0^T\hspace{-5pt}\dd t\,\ee_\nu\,[\int_0^{t}\hspace{-1pt}
V(s,\eta_{cs})\,\dd s\,]^2
\,\le\,
\frac{2e^{\beta T}}{\beta c}\,
\hspace{-3pt}\int_0^T\hspace{-5pt}
\left(\rule{0pt}{10pt}\right.\!
V(s,\cdot)
\;\rule[-4pt]{1pt}{14pt}\;
(\frac{\beta}{2c}-L^{\rm sym})^{-1}V(s,\cdot)
\!\left.\rule{0pt}{10pt}\right)_{L^2(\nu)}\dd s$$
for all bounded measurable functions $V$ 
on $[0,\infty)\times X$.
\end{itemize}
\end{lemma}
{\it Proof}. The conditions given in (i) are obvious conditions
for $L_s\,g(s,\eta)=[L^{\rm sym}g(s,\cdot)](\eta)$
to be true which indeed proves the claim.

Part (ii) only needs to be proven for 
$V(s,\eta){\bf 1}_{[0,T]}(s)$ instead of $V(s,\eta)$.
The method is to approximate
$V{\bf 1}_{[0,T]}$ by `good' functions $V_n$ 
such that the functions $g_n$ corresponding
to $\tilde{V}_n$ satisfy the conditions of part (i).
Going through the proof of Corollary \ref{my inequ} 
but using the conditions of part (i) reveals that
one {\em can} conclude that
\begin{equation}\label{one can}
(\tilde{V}_n\,|\,
G_{\!\frac{\beta}{2c}}\tilde{V}_n)_{L^2(\ell\otimes\nu)}
\,\le\,
(\tilde{V}_n\,|\,
(\frac{\beta}{2c}-L_s)^{-1}\tilde{V}_n)_{L^2(\ell\otimes\nu)}
\end{equation}
in this specific case
where, by part (i), the right-hand side is equal to
\begin{eqnarray*}
&&\int_0^\infty\hspace{-1pt}
(\tilde{V}_n(s,\cdot)\,|\,
(\frac{\beta}{2c}-L^{\rm sym})^{-1}\tilde{V}_n(s,\cdot))_{L^2(\nu)}\,\dd s\\
&=&\rule{0pt}{13pt}
c\int_0^\infty\hspace{-1pt}e^{-\beta s}\,
({V}_n(s,\cdot)\,|\,
(\frac{\beta}{2c}-L^{\rm sym})^{-1}{V}_n(s,\cdot))_{L^2(\nu)}\,\dd s
\end{eqnarray*}
so that
\begin{equation}\label{stand approx}
\int_0^T\hspace{-5pt}\dd t\,\ee_\nu\,[\int_0^{t}\hspace{-1pt}
V_n(s,\eta_{cs})\,\dd s\,]^2
\,\le\,
\frac{2e^{\beta T}}{\beta c}\,
\hspace{-3pt}\int_0^\infty\hspace{-5pt}
\left(\rule{0pt}{10pt}\right.\!
V_n(s,\cdot)
\;\rule[-4pt]{1pt}{14pt}\;
(\frac{\beta}{2c}-L^{\rm sym})^{-1}V_n(s,\cdot)
\!\left.\rule{0pt}{10pt}\right)_{L^2(\nu)}\dd s
\end{equation}
for all $n$ by Lemma \ref{lem_inhomo}. Taking limits when $n$ goes to
infinity in the above inequality finally proves part (ii).
Notice that the upper limit of the $\dd s$-integration 
can indeed be changed to $T$ 
because $V_n$ approximates $V{\bf 1}_{[0,T]}$.
\hfill$\Box$
\begin{rem}\rm\label{smooth}
\begin{itemize}\item[(i)]
Taking limits in
(\ref{stand approx}) makes clear that the approximation of
$V{\bf 1}_{[0,T]}$ by `good' functions $V_n$ is a standard
approximation of a function in $L^2(\ell\otimes\nu)$ by 
in some sense `smooth' functions and the details are therefore omitted.
\item[(ii)]
Remember that the left-hand side of (\ref{one can}) can be transformed into 
$$\int_0^{\infty}\hspace{-5pt}\dd s\int_0^{\infty}\hspace{-5pt}\dd r
\,e^{-\beta(s+r)/c}\,
\ee_\nu\,{V}_n(s/c,\eta_{s}){V}_n((s+r)/c,\eta_{s+r})$$
as in the proof of Lemma \ref{lem_inhomo}.
However, it is not possible to transform the right-hand side of
(\ref{one can}) into a similar expression because the resolvent
$((\alpha-L_s)^{-1})_{\alpha>0}$ is not associated with a process of
the form $(s,\eta^{\rm sym}_s)_{s\ge 0}$.
\item[(iii)]
In the case were $V$ is time independent one can directly apply 
Corollary \ref{my inequ} to the right-hand side of the inequality in
Remark \ref{no t remark}(i) which gives
$$\int_0^T\hspace{-5pt}\dd t\,\ee_\nu\,[\int_0^{t}\hspace{-1pt}
V(\eta_{cs})\,\dd s\,]^2
\,\le\,
\frac{2e^{\beta T}}{\beta^2 c}\,
({V}\,|\,
(\beta/c-L^{\rm sym})^{-1}{V})_{L^2(\nu)}$$
for all bounded $V$ in $L^2(\nu)$ hence for all 
$V\in L^2(\nu)$ by approximation.
Notice  that $\nu$ is an invariant measure for
the resolvent $((\alpha-L^{\rm sym})^{-1})_{\alpha>0}$ 
under the conditions made.
\end{itemize}
\end{rem}
\section{Application to 1-dimensional
Simple Exclusion}\label{section application}
Fix $p,q\ge 0$ such that $p+q=1$ and
let $(\Omega,{\cal F},\pp_{\!\!\eta},\eta\in\{0,1\}^\NZ,
(\eta_t)_{t\ge 0})$ denote the strong Markov Feller process
whose generator $L$
acts on local functions $f:\{0,1\}^\NZ\to\NR$ as
{\small
\begin{equation}\label{operator}
Lf(\eta)\,=
\sum_{x\in\NZ}\left(
2p\,\eta(x)(1-\eta(x+1))[f(\eta^{x,x+1})-f(\eta)]
+\;2q\,\eta(x)(1-\eta(x-1))[f(\eta^{x,x-1})-f(\eta)]
\right)
\end{equation}
}

\noindent
where the operation
$$\eta^{x,y}(z)\,=\,\left\{\begin{array}{ccc}
\eta(z)&:&z\not=x,y\\
\eta(x)&:&z=y\\
\eta(y)&:&z=x
\end{array} \right.$$
exchanges the ``spins'' at $x$ and $y$. 
This process is called simple exclusion process, see \cite{L1999}
for a good account on the existing theory.

Denote by $\nu_{1/2}$ the Bernoulli product
measure on $\{0,1\}^{\NZ}$ satisfying
$\nu_{1/2}(\eta(x)=1)=1/2$ for all $x\in\NZ$ which is one of the
invariant ergodic states of the simple exclusion process. If
$$\pp\,=\,\int\pp_{\!\!\eta}\,\dd\nu_{1/2}(\eta)
\quad\mbox{as well as}\quad
\xi_t(x)\,=\,\frac{\eta_t(x)-\ee\eta_t(x)}
{\sqrt{\mbox{\bf Var}(\eta_t(x))}}$$
where $\ee$ and $\mbox{\bf Var}$ stand for expectation and
variance with respect to \pp, respectiveley,
then the process $(\xi_t)_{t\ge 0}$ is a stationary process on 
$(\Omega,{\cal F},\pp)$ which takes values in $\{-1,1\}^{\NZ}$
and the push forward of $\nu_{1/2}$ with respect to the map
$$\eta\mapsto\xi
\quad\mbox{given by}\quad
\xi(x)\,=\,\frac{\eta(x)-1/2}{\sqrt{1/4}},\;x\in\NZ,$$
is the invariant distribution of $\xi_t,\,t\ge 0$.

For $\Lambda\subseteq\NZ$ finite, set 
$\xi_\Lambda=\prod_{x\in\Lambda}\xi(x)$ if $\Lambda$ is not empty and
$\xi_\emptyset=1$ otherwise. Then,
$\{\xi_\Lambda:\mbox{$\Lambda\subseteq\NZ$ finite}\}$ forms an
orthonormal basis of $L^2(\nu_{1/2})$. 
Hence the linear hull $\LL\{\xi_\Lambda\}$
of $\{\xi_\Lambda:\mbox{$\Lambda\subseteq\NZ$ finite}\}$
is dense in $L^2(\nu_{1/2})$.
Remark that the operator $(L,\LL\{\xi_\Lambda\})$
is closable on $L^2(\nu_{1/2})$
and that its closure $(L,D(L))$
generates a Markovian strongly continuous contraction semigroup 
$(T_t)_{t\ge 0}$ on $L^2(\nu_{1/2})$ which is associated with the
transition semigroup of the strong Markov process 
$(\Omega,{\cal F},\pp_{\!\!\eta},\eta\in\{0,1\}^\NZ,
(\eta_t)_{t\ge 0})$. 
Furthermore, it follows from (\ref{operator}) that
$$L\xi_\Lambda\,=\,
\gamma{\cal A}_+\xi_\Lambda-\gamma{\cal A}_+^\ast\xi_\Lambda
+{\cal S}\xi_\Lambda
\quad\mbox{where}\quad
\gamma=p-q$$
and
$$  
{\cal A}_+\xi_\Lambda\,=\,
\sum_{x\in\Lambda}\left[
{\bf 1}_{\Lambda^c}(x+1)\xi_{\Lambda\cup\{x+1\}}
-{\bf 1}_{\Lambda^c}(x-1)\xi_{\Lambda\cup\{x-1\}}
\right].
$$  
The adjoint operator of $({\cal A}_+,\LL\{\xi_\Lambda\})$ with respect
to the inner product on $L^2(\nu_{1/2})$ is denoted by ${\cal A}_+^\ast$.
Its domain includes $\LL\{\xi_\Lambda\}$ and, on this subdomain, it is
given by
$$  
{\cal A}_+^\ast\xi_\Lambda\,=\,
\sum_{x\in\Lambda}\left[
{\bf 1}_{\Lambda^c}(x+1)-{\bf 1}_{\Lambda^c}(x-1)
\right]\xi_{\Lambda\setminus\{x\}}\,.
$$  
The operator $({\cal S},\LL\{\xi_\Lambda\})$ is a symmetric operator
on $L^2(\nu_{1/2})$ satisfying
$$\hspace{-7.7cm}
{\cal S}\xi_\Lambda\,=\,
{\cal S}_0\xi_\Lambda-2|\Lambda|\xi_\Lambda$$
where 
\begin{eqnarray*}
{\cal S}_0\xi_\Lambda&=&
\sum_{x\in\Lambda}\left[
{\bf 1}_{\Lambda^c}(x+1)\xi_{\Lambda\setminus\{x\}\cup\{x+1\}}
+{\bf 1}_{\Lambda}(x+1)\xi_{\Lambda}\right.\\
&&\hspace{3cm}\left.
+{\bf 1}_{\Lambda^c}(x-1)\xi_{\Lambda\setminus\{x\}\cup\{x-1\}}
+{\bf 1}_{\Lambda}(x-1)\xi_{\Lambda}
\right]
\end{eqnarray*}
and $|\Lambda|$ denotes the cardinality of $\Lambda$.

As a consequence, the adjoint $(L^\ast,D(L^\ast))$ of $(L,D(L))$
satisfies
$$L^\ast\,=\,\gamma{\cal A}_+^\ast-\gamma{\cal A}_++{\cal S}
\quad\mbox{on}\quad
\LL\{\xi_\Lambda\}$$
and $(L^\ast,D(L^\ast))$ is the closure of
$(L^\ast,\LL\{\xi_\Lambda\})$. Thus if
$$G_\alpha\,=\,\int_0^\infty e^{-\alpha t}T_t\,\dd t,
\quad\alpha>0,$$
is the Markovian strongly continuous contraction resolvent generated
by $(L,D(L))$ then $(L^\ast,D(L^\ast))$ generates $G_\alpha^\ast$
hence the condition $(A1)$ in Section \ref{section resolv method} is
satisfied. Furthermore
$$L_s\,=\,\frac{L+L^\ast}{2}\,=\,S
\quad\mbox{on}\quad
\LL\{\xi_\Lambda\}$$
and $(S,\LL\{\xi_\Lambda\})$ is of course essentially self-adjoint,
its closure being the generator of the symmetric simple exclusion process.
Altogether the domain $\LL\{\xi_\Lambda\}$ satisfies 
all the conditions imposed on the domain $D_0$ in Corollary \ref{my inequ}.
\begin{rem}\rm\label{not sym part}
The operator $(\gamma{\cal A}_++{\cal S},\LL\{\xi_\Lambda\})$
is an easy example of an operator which, when seen as a perturbation of 
$({\cal S},\LL\{\xi_\Lambda\})$
does not have $S$ as its symmetric part on $\LL\{\xi_\Lambda\}$.
\end{rem}

The example field of type (\ref{field}) discussed in this paper
is the density fluctuation field
$$Y_s^\varepsilon\,=\,\sqrt{\varepsilon}
\sum_{x\in\NZ}\xi_{s\varepsilon^{-2}}(x)
\delta_{\varepsilon x},\quad s\ge 0,$$
in the case of $\sqrt{\varepsilon}$-asymmetric one-dimensional simple 
exclusion where $p$ and $q$ hence $\gamma$ depend on $\varepsilon$
such that
$$\gamma\,=\,\gamma_\varepsilon\,=\,\tilde{\gamma}\cdot\sqrt{\varepsilon}.$$
As a consequence $L,T_t,G_\alpha$ and $\ee$ introduced above
are denoted by 
$L_\varepsilon,T_t^\varepsilon,G_\alpha^\varepsilon$ and
$\ee_\varepsilon$ in what follows.

If $\varepsilon$ is small then the field
has a related approximate equation of type (\ref{approx equ})
with ${\cal A}$ being the one-dimensional Laplacian and with
$V_\varepsilon^G(s,\cdot)$
being of the special form $\tilde{\gamma}V_\varepsilon^G$ where
$$V_\varepsilon^G(\xi)\,=\,
-\sum_{x\in\NZ}G^\prime(\varepsilon x)\xi(x)\xi(x+1).$$
Functions of this type are called quadratic fluctuations in this
paper, {\it quadratic} since its summands are of type $\xi_\Lambda$
with $|\Lambda|=2$ and {\it fluctuations} since $\xi$ is
$(\eta-1/2)/\sqrt{1/4}$. Remark that $V_\varepsilon^G$
does not depend on $s$ because of the special choice of $\nu_{1/2}$
to be the invariant state of the system. However, the choice of a
different invariant state would only make the notation more
complicated but would not affect the arguments used below.

Furthermore, $V_\varepsilon^G$ remains bounded if $G$ is a test
function in the Schwartz space ${\mathscr S}(\NR)$ and this space of
test functions is chosen in what follows.

The functional $F$ used to replace
$$\int_0^t V_\varepsilon^G(\xi_{s\varepsilon^{-2}})\,\dd s
\quad\mbox{by}\quad
\int_0^t F(Y_s^\varepsilon,G)\,\dd s$$
depends on a further parameter $N$, that is $F=F_N$ in this example.
It can be defined as follows.
Fix an \underline{even}
non-negative test function $d$ satisfying
$$\mbox{supp}\,d=[-1,1]\quad\mbox{and}\quad
\int_\NR d(u)\,\dd u=1$$
and denote by $d_N$ the function $x\mapsto Nd(Nx),\,N\ge 1$.
Remark that the convolution $Y^{\varepsilon}_s\star d_N$ is a
$C^\infty$-function on $\NR$ satisfying
$$\int_\NR H(u)\,\dd Y^{\varepsilon}_s(u)
\,=\,\lim_{N\uparrow\infty}
\int_\NR H(u)(Y^{\varepsilon}_s\star d_N)(u)\,\dd u$$
for all $s\ge 0$ and all continuous functions $H$ on $\NR$
with sufficiently fast decaying tails.
Now set
$$F_N({\cal Y},G)\,=\,
-\int_\NR G^\prime(u)({\cal Y}\star d_N)^2(u)\,\dd u
,\quad{\cal Y}\in{\mathscr S}^\prime(\NR).$$
Fixing a finite time horizon $T>0$,
one wants to estimate
$$\int_0^T\hspace{-3pt}\dd t\,
\ee_\varepsilon\left(
\int_0^t F_N(Y_s^\varepsilon,G)\,\dd s
\,-\,
\int_0^t V_\varepsilon^G(\xi_{s\varepsilon^{-2}})\,\dd s
\right)^{\hspace{-2pt}2}$$
for small $\varepsilon$ and large $N$.
So fix $\varepsilon,N$ and observe that
$$\int_0^t F_N(Y_s^\varepsilon,G)\,\dd s
\,-\,
\int_0^t V_\varepsilon^G(\xi_{s\varepsilon^{-2}})\,\dd s$$
$$=-\int_0^t\!\int_\NR
G^\prime(u)(Y^{\varepsilon}_s\star d_N)^2(u)\,\dd u\dd s
\,+
\int_0^t
\sum_{x\in\NZ}G^\prime(\varepsilon x)
\xi_{s\varepsilon^{-2}}(x)\xi_{s\varepsilon^{-2}}(x+1)
\,\dd s$$
hence, setting $H=-G^\prime$, one can split
$$\int_0^t F_N(Y_s^\varepsilon,G)\,\dd s
\,-\,
\int_0^t V_\varepsilon^G(\xi_{s\varepsilon^{-2}})\,\dd s
\quad\mbox{into}\quad
\sum_{i=1}^4
\int_0^t
V_{\varepsilon,N}^{H,i}(\xi_{s\varepsilon^{-2}})\,\dd s$$
using further quadratic fluctuations
$V_{\varepsilon,N}^{H,1},V_{\varepsilon,N}^{H,2},V_{\varepsilon,N}^{H,3},
V_{\varepsilon,N}^{H,4}$
given by
\begin{eqnarray*}
V_{\varepsilon,N}^{H,1}(\xi)&=&
\sum_{x\in\NZ}
\int_\NR
[H(u)-H(\varepsilon x)]
d_N(u-\varepsilon x)
\sum_{\tilde{x}\in\NZ}
\varepsilon
d_N(u-\varepsilon\tilde{x})\,\dd u\;
\xi(x)\xi(\tilde{x}),\\
V_{\varepsilon,N}^{H,2}(\xi)&=&
\varepsilon
\sum_{x\in\NZ}H(\varepsilon x)
\int_\NR d_N^2(u-\varepsilon x)\,\dd u\;
\xi(x)[\xi(x)-\xi(x+1)],\\
V_{\varepsilon,N}^{H,3}(\xi)&=&
\varepsilon
\sum_{x\not=\tilde{x}}H(\varepsilon x)
\int_\NR d_N(u-\varepsilon x)d_N(u-\varepsilon\tilde{x})\,\dd u\;
\xi(x)[\xi(\tilde{x})-\xi(x+1)],\\
V_{\varepsilon,N}^{H,4}(\xi)&=&
\sum_{x\in\NZ}H(\varepsilon x)
\int_\NR d_N(u-\varepsilon x)
\hspace{-3pt}\left[
\sum_{\tilde{x}\in\NZ}\varepsilon d_N(u-\varepsilon\tilde{x})-1
\right]\!\dd u\;
\xi(x)\xi(x+1).
\end{eqnarray*}
Hence, if $(\cdot\,|\,\cdot)$ denotes the inner product on
$L^2(\nu_{1/2})$ then
\begin{eqnarray}
&&\ee_\varepsilon\left(
\int_0^t F_N(Y_s^\varepsilon,G)\,\dd s
\,-\,
\int_0^{t}
V_\varepsilon^G(\xi_{s\varepsilon^{-2}})\,\dd s
\right)^{\hspace{-2pt}2}\nonumber\\
&\le&\rule{0pt}{30pt}
4\sum_{i=1}^4
\ee_\varepsilon\left(
\int_0^t
V_{\varepsilon,N}^{H,i}(\xi_{s\varepsilon^{-2}})\,\dd s
\right)^{\hspace{-2pt}2}\label{four summands}\\
&\le&\rule{0pt}{30pt}
\label{split difference}
8\sum_{i=1}^3\,
\varepsilon^4\hspace{-4pt}
\int_0^{t\varepsilon^{-2}}\hspace{-12pt}\dd s
\int_0^s\dd r\,
(V_{\varepsilon,N}^{H,i}\,|\,T_r^\varepsilon V_{\varepsilon,N}^{H,i})
\,+\,
4\,\ee_\varepsilon\left(\varepsilon^2\hspace{-4pt}
\int_0^{t\varepsilon^{-2}}\hspace{-12pt}
V_{\varepsilon,N}^{H,4}(\xi_s)\,\dd s
\right)^{\hspace{-2pt}2}
\end{eqnarray}
for all $t\ge 0$ by the Markov property.

In what follows, $\|H\|_p$
denotes the norm of a test function $H$ in $L^p(\NR),\,1\le p\le\infty$.
Then, as
$|\sum_{\tilde{x}\in\NZ}\varepsilon d_N(u-\varepsilon\tilde{x})-1|
\le\varepsilon^2 N^2\|d^{\prime\prime}\|_\infty$ is easily realised,
one obtains that
\begin{equation}\label{estiV4}
\ee_\varepsilon\left(\varepsilon^2\hspace{-4pt}
\int_0^{t\varepsilon^{-2}}\hspace{-12pt}
V_{\varepsilon,N}^{H,4}(\xi_s)\,\dd s
\right)^{\hspace{-2pt}2}
\le\,
\varepsilon^2 N^4\|d^{\prime\prime}\|_\infty^2\|H\|_1^2\cdot t^2,
\quad t\ge 0.
\end{equation}
The other three integrals in (\ref{split difference}) 
are much harder to control. But the following lemma gives bounds for
these integrals if $T_r^\varepsilon$ is substituted by the semigroup
$T_r^{\rm sym}$ associated with the symmetric simple exclusion
process.
\begin{lemma}\label{esti_sym}
Let $H\in{\mathscr S}(\NR)$ be a test function such that
$\int_\NR H(u)\,\dd u=0$. Then 
{\small
$$\begin{array}{crcl}
(i)\hspace{.2cm}&\displaystyle
\varepsilon^4\hspace{-4pt}
\int_0^{t\varepsilon^{-2}}\hspace{-12pt}\dd s
\int_0^s\dd r\,
(V_{\varepsilon,N}^{H,1}\,|\,T_r^{\rm sym}V_{\varepsilon,N}^{H,1})
&\le&\displaystyle
c_0\,t^2
\left(\rule{0pt}{11pt}
\|d\|_2^4\,\frac{\|(1+u^2)H^{\prime\prime}\|_\infty^2}{N^2}
+\|d\|_\infty^2\,\frac{\|(1+u^2)H^\prime\|_\infty^2}{N}
\right)\\
\rule{0pt}{30pt}
(ii)\hspace{.2cm}&\displaystyle
\varepsilon^4\hspace{-4pt}
\int_0^{t\varepsilon^{-2}}\hspace{-12pt}\dd s
\int_0^s\dd r\,
(V_{\varepsilon,N}^{H,2}\,|\,T_r^{\rm sym}V_{\varepsilon,N}^{H,2})
&\le&\displaystyle
t^2\,\|d\|_2^4\,
\left(\rule{0pt}{11pt}
\varepsilon^2 N^2\,\|H^\prime\|_\infty^2
+\varepsilon N^2\,\|(1+u^2)H\|_\infty^2
\right)\\
\rule{0pt}{30pt}
(iii)\hspace{.2cm}&\displaystyle
\varepsilon^4\hspace{-4pt}
\int_0^{t\varepsilon^{-2}}\hspace{-12pt}\dd s
\int_0^s\dd r\,
(V_{\varepsilon,N}^{H,3}\,|\,T_r^{\rm sym}V_{\varepsilon,N}^{H,3})
&\le&\displaystyle
c_0\,t^2\,\|d\|_\infty^2\,\frac{\|(1+u^2)H\|_\infty^2}{N^{1/3}}
\end{array}$$
}

\noindent
for all $t\ge 0,\,N=1,2,\dots,\,\varepsilon>0$
where $c_0$ is a constant which neither depends 
on the chosen test function $H$ nor on the mollifier $d$.
\end{lemma}
{\it Proof}. For showing $(i)$ one splits $V_{\varepsilon,N}^{H,1}$
into two sums
\begin{eqnarray}\label{splitV1}
&&\varepsilon
\sum_{x\in\NZ}
\int_\NR
[H(u)-H(\varepsilon x)]
d_N^2(u-\varepsilon x)\,\dd u\nonumber\\
&+&
\sum_{x\not=\tilde{x}}
\int_\NR
[H(u)-H(\varepsilon x)]
d_N(u-\varepsilon x)\,
\varepsilon d_N(u-\varepsilon\tilde{x})\,\dd u\;
\xi(x)\xi(\tilde{x}).
\end{eqnarray}
Applying the Taylor expansion
$$H(u)-H(\varepsilon x)\,=\,(u-\varepsilon x)H^\prime(\varepsilon x)
\,+\,(u-\varepsilon x)^2 H^{\prime\prime}(\theta_{\varepsilon x}^u)/2
\quad\mbox{with}\quad
\theta_{\varepsilon x}^u\in[\varepsilon x-u,\varepsilon x+u]$$
to the first sum yields
$$
\left|
\varepsilon
\sum_{x\in\NZ}
\int_\NR
[H(u)-H(\varepsilon x)]
d_N^2(u-\varepsilon x)\,\dd u
\right|
\,\le\,
\varepsilon\hspace{-4pt}
\sum_{x\in\NZ}\int_\NR
\frac{(u-\varepsilon x)^2}{2}
|H^{\prime\prime}(\theta_{\varepsilon x}^u)|
d_N^2(u-\varepsilon x)\,\dd u\\
$$
where $\int_\NR(u-\varepsilon x)d_N^2(u-\varepsilon x)\dd u$
vanishes because the mollifier $d$ is even. 
Now observe that
$$|H^{\prime\prime}(\theta_{\varepsilon x}^u)|
\,\le\left\{\begin{array}{rcl}
\sup_{\tilde{u}}|(1+\tilde{u}^2)H^{\prime\prime}(\tilde{u})|
\cdot[1+(\varepsilon x-\frac{1}{N})^2]^{-1}&:&x>1/\varepsilon\\
\rule{0pt}{20pt}
\sup_{\tilde{u}}|(1+\tilde{u}^2)H^{\prime\prime}(\tilde{u})|
\cdot[1+(\varepsilon x+\frac{1}{N})^2]^{-1}&:&x<-1/\varepsilon
\end{array}\right.$$
on the set $\{u\in\NR:d_N(u-\varepsilon x)\not=0\}$
since ${\rm supp}\,d_N=[-\frac{1}{N},\frac{1}{N}]$.
Furthermore
$$(u-\varepsilon x)^2\,\le\,1/N^2
\quad\mbox{on}\quad\{u\in\NR:d_N(u-\varepsilon x)\not=0\}$$
and
$$\int_\NR d_N^2(u-\varepsilon x)\,\dd u\,=\,N\|d\|_2^2
\quad\mbox{for all $x\in\NZ$.}$$
Therefore, the sum immediately above (\ref{splitV1})
can be estimated by
{\small
$$\frac{\|d\|_2^2}{2N}\left(\rule{0pt}{13pt}\right.
\|(1+{u}^2)H^{\prime\prime}\|_\infty
\overbrace{
\sum_{x<-\frac{1}{\varepsilon}}\varepsilon
[1+(\varepsilon x+\frac{1}{N})^2]^{-1}}^{\le \pi/2}
+\;2\|H^{\prime\prime}\|_\infty\,+\;
\|(1+{u}^2)H^{\prime\prime}\|_\infty
\overbrace{
\sum_{x>\frac{1}{\varepsilon}}\varepsilon
[1+(\varepsilon x-\frac{1}{N})^2]^{-1}}^{\le \pi/2}
\left.\rule{0pt}{13pt}\right)$$
}
\begin{equation}\label{with pi}
\le\,
\frac{\|d\|_2^2}{2N}\,(\pi+2)\|(1+{u}^2)H^{\prime\prime}\|_\infty
\end{equation}
which explains the first summand on the right-hand side of $(i)$. 

The second summand is a bound of the integral on the left-hand side
of $(i)$ but with $V_{\varepsilon,N}^{H,1}$ replaced by the
fluctuation field given by (\ref{splitV1}). This bound 
is obtained by copying the proof of Lemma 1 in \cite{A2007} for
$x_0=1$ using the equality
$$H(u)-H(\varepsilon x)\,=\,(u-\varepsilon x)
H^\prime(\tilde{\theta}_{\varepsilon x}^u)$$
where
$\tilde{\theta}_{\varepsilon x}^u\in[\varepsilon x-u,\varepsilon x+u]$.
The only difference to the proof of Lemma 1 in \cite{A2007} is that,
similar to how (\ref{with pi}) was derived, the sum
$\sum_{x\in\NZ}\varepsilon|H^\prime(\tilde{\theta}_{\varepsilon x}^u)|$
is estimated by 
$(\pi+2)\sup_{\tilde{u}}|(1+\tilde{u}^2)H^{\prime}(\tilde{u})|$
and not by $c_H\|H^\prime\|_\infty$.

The left-hand side of (ii) can be estimated 
the same way the sum $S_1(t,\varepsilon,N)$ was estimated 
in the proof of Theorem 1 in \cite{A2007}. 
Following this proof would give
$$\varepsilon^4\hspace{-4pt}
\int_0^{t\varepsilon^{-2}}\hspace{-12pt}\dd s
\int_0^s\dd r\,
(V_{\varepsilon,N}^{H,2}\,|\,T_r^{\rm sym}V_{\varepsilon,N}^{H,2})
\,\le\,
t^2\,\|d\|_2^4\,
\left(\rule{0pt}{11pt}
N^2
(\sum_{x\in\NZ}\varepsilon H(\varepsilon x))^2
+\varepsilon N^2\,c_H\|H\|_\infty^2
\right)$$
if $H$ were a test function with compact support.
But this implies $(ii)$ because,
by our assumption $\int_\NR H(u)\,\dd u=0$,
it holds that
$|\sum_{x\in\NZ}\varepsilon H(\varepsilon x)|
\le\varepsilon\|H^\prime\|_\infty$ by our assumption
$\int_\NR H(u)\,\dd u=0$. Again, as in the proof of part $(i)$
above, the bound $c_H\|H\|_\infty$ is replaced by
$\|(1+u^2)H\|_\infty$.

Finally, the inequality $(iii)$ can be established
by copying the proof of Theorem 1 in \cite{A2007}
with respect to $S_{6-9}(t,\varepsilon,N)$ for $x_0=1$ and $\alpha=2/3$ 
manipulating the constant $c_H$ accordingly.\hfill$\Box$\\

The next lemma is the key result of this section. It translates the
estimates given by Lemma \ref{esti_sym} into estimates of the summands
in (\ref{four summands}) on page \pageref{four summands}
when integrating them against $\dd t$ over $t\in[0,T]$.
\begin{lemma}\label{key result}
Let $H\in{\mathscr S}(\NR)$ be a test function such that
$\int_\NR H(u)\,\dd u=0$ and fix a finite time horizon $T>0$. Then
$$\begin{array}{crcl}
(i)\hspace{.2cm}&\displaystyle
\int_0^T\hspace{-3pt}\dd t\,
\ee_\varepsilon\left(
\int_0^t
V_{\varepsilon,N}^{H,1}(\xi_{s\varepsilon^{-2}})\,\dd s
\right)^{\hspace{-2pt}2}
&\le&\displaystyle
e^T C_d
\left(\rule{0pt}{11pt}
\frac{\|(1+u^2)H^{\prime\prime}\|_\infty^2}{N^2}
+\frac{\|(1+u^2)H^\prime\|_\infty^2}{N}
\right)\\
\rule{0pt}{30pt}
(ii)\hspace{.2cm}&\displaystyle
\int_0^T\hspace{-3pt}\dd t\,
\ee_\varepsilon\left(
\int_0^t
V_{\varepsilon,N}^{H,2}(\xi_{s\varepsilon^{-2}})\,\dd s
\right)^{\hspace{-2pt}2}
&\le&\displaystyle
e^T C_d
\left(\rule{0pt}{11pt}
\varepsilon^2 N^2\,\|H^\prime\|_\infty^2
+\varepsilon N^2\,\|(1+u^2)H\|_\infty^2
\right)\\
\rule{0pt}{30pt}
(iii)\hspace{.2cm}&\displaystyle
\int_0^T\hspace{-3pt}\dd t\,
\ee_\varepsilon\left(
\int_0^t
V_{\varepsilon,N}^{H,3}(\xi_{s\varepsilon^{-2}})\,\dd s
\right)^{\hspace{-2pt}2}
&\le&\displaystyle
e^T C_d
\,\frac{\|(1+u^2)H\|_\infty^2}{N^{1/3}}\\
\rule{0pt}{30pt}
(iv)\hspace{.2cm}&\displaystyle
\int_0^T\hspace{-3pt}\dd t\,
\ee_\varepsilon\left(
\int_0^t
V_{\varepsilon,N}^{H,4}(\xi_{s\varepsilon^{-2}})\,\dd s
\right)^{\hspace{-2pt}2}
&\le&\displaystyle
T^3 C_d
\,\varepsilon^2 N^4\|H\|_1^2
\end{array}$$
for all $N=1,2,\dots,\,\varepsilon>0$ where $C_d$ is a constant
which only depends on the choice of the mollifier $d$.
\end{lemma}
{\it Proof}.
Choose $\beta=1$ and consider $i=1,2,3$ first.
Then, applying the inequality in Remark \ref{smooth}(iii)
with respect to $c=\varepsilon^{-2}$, one obtains that
$$\int_0^T\hspace{-3pt}\dd t\,
\ee_\varepsilon\left(
\int_0^t
V_{\varepsilon,N}^{H,i}(\xi_{s\varepsilon^{-2}})\,\dd s
\right)^{\hspace{-2pt}2}
\,\le\,
2e^{T}\varepsilon^2\,
(V_{\varepsilon,N}^{H,i}\,|\,
(\varepsilon^2-L^{\rm sym})^{-1}V_{\varepsilon,N}^{H,i})$$
where
$$(V_{\varepsilon,N}^{H,i}\,|\,
(\varepsilon^2-L^{\rm sym})^{-1}V_{\varepsilon,N}^{H,i})
\,=\,
\varepsilon^2
\int_0^\infty\hspace{-3pt}\dd t\,e^{-t}
\int_0^{t\varepsilon^{-2}}\hspace{-12pt}\dd s
\int_0^s\dd r\,
(V_{\varepsilon,N}^{H,i}\,|\,T_r^{\rm sym}V_{\varepsilon,N}^{H,i})$$
by a calculation similar to the proof of Lemma \ref{lem_inhomo}
in the time independent case.
Hence (i)-(iii) follows from Lemma \ref{esti_sym} since
$\int_0^\infty\hspace{-1pt}t^2e^{-t}\,\dd t$ is finite.

Finally, (iv) follows directly from (\ref{estiV4}) by integration against
$\dd t$ over $t\in[0,T]$.\hfill$\Box$\\

Applying Lemma \ref{key result} in the case where $H$ is taken to be
$-G^\prime$ immediately gives the corollary below. Notice that 
$\int_\NR H(u)\,\dd u=0$ is of course satisfied for $H=-G^\prime$.
\begin{cor}\label{weak replacement}
Fix an arbitrary but finite time horizon $T>0$. 
Then, for every smooth test function $G\in{\mathscr S}(\NR)$, it holds that
$$\lim_{\rule{0pt}{8pt}N\uparrow\infty}\;
\limsup_{\varepsilon\downarrow 0}
\int_0^T\hspace{-3pt}\dd t\,
\ee_\varepsilon\left(
\int_0^t F_N(Y_s^\varepsilon,G)\,\dd s
\,-\,
\int_0^t V_\varepsilon^G(\xi_{s\varepsilon^{-2}})\,\dd s
\right)^{\hspace{-2pt}2}
=\,0$$ 
where
$$F_N({\cal Y},G)\,=\,
-\int_\NR G^\prime(u)({\cal Y}\star d_N)^2(u)\,\dd u
\quad\mbox{and}\quad
V_\varepsilon^G(\xi)\,=\,
-\sum_{x\in\NZ}G^\prime(\varepsilon x)\xi(x)\xi(x+1)$$
for ${\cal Y}\in{\mathscr S}^\prime(\NR)$ and $\xi\in\{-1,1\}^\NZ$, respectively.
\end{cor}
\begin{rem}\rm
A replacement of 
$$\int_0^t V_\varepsilon^G(\xi_{s\varepsilon^{-2}})\,\dd s
\quad\mbox{by}\quad
\int_0^t F_N(Y_s^\varepsilon,G)\,\dd s$$
for every $t\in[0,T]$ and not only for an average over $t\in[0,T]$ was shown
in \cite{JG2010} for a slightly different functional $F_N$. However it
has been shown in \cite{A2012} that many of the conclusions drawn 
from the stronger 
replacement result\footnote{This replacement result was called
`Second-order Boltzmann-Gibbs principle'.}
in \cite{JG2010} can actually be
obtained by only applying the weaker replacement result
of Corollary \ref{weak replacement}.
\end{rem}

\end{document}